\documentclass[12pt]{article}

\usepackage{amssymb,hyperref}

\def\wPms{{$\widehat P$-matrices}}

\begin{document}
\newcommand\Z{{\mathbb Z}}
\newcommand\R{{\mathbb R}}




\begin{center}
{\large Canonical bases of invariant polynomials for the irreducible reflection groups of types $E_6$, $E_7$, and $E_8$}\\
\vskip 0.5cm
{\large Vittorino Talamini}\\
{\tt vittorino.talamini@uniud.it}\\

\end{center}



\begin{abstract}
Given a rank $n$ irreducible finite reflection group $W$, the $W$-invariant polynomial functions defined in $\R^n$ can be written as polynomials of $n$ algebraically independent homogeneous polynomial functions, $p_1(x),\ldots,p_n(x)$, called basic invariant polynomials. Their degrees are well known and typical of the given group $W$. The polynomial $p_1(x)$ has the lowest degree, equal to 2. It has been proved that it is possible to choose all the other $n-1$ basic invariant polynomials in such a way that they satisfy a certain system of differential equations, including the Laplace equations $\triangle p_a(x)=0,\ a=2,\ldots,n$, and so are harmonic functions. Bases of this kind are called canonical. Explicit formulas for canonical bases of invariant polynomials have been found for all irreducible finite reflection groups, except for those of types $E_6$, $E_7$ and $E_8$. Those for the groups of types $E_6$, $E_7$ and $E_8$ are determined in this article.
\end{abstract}








\section{Introduction}

When a rank $n$ irreducible finite reflection group $W$ acts in $\R^n$, there are no fixed points besides the origin of $\R^n$, and, with no loss of generality, one can assume that  $W\subset O(n)$. By saying that $f(x)$, $x\in\R^n$ is a {\em $W$-invariant function}, one means that $f(gx)=f(x)$, for all  $g\in W$, $x\in\R^n$. It is well known from Ref.~\cite{Chev1955
} that there exist a basis of $n$ algebraically independent $W$-invariant real homogeneous polynomial functions $p_1(x),\ldots,p_n(x)$, called {\em basic invariant polynomials}, such that all $W$-invariant polynomial 
functions can be univocally written as polynomial 
functions of the basic invariant polynomials: $f(x)=\widehat f(p_1(x),\ldots,p_n(x))$,  with $\widehat f(p_1,\ldots,p_n)=\widehat f(p)$ a polynomial 
function of $p\in\R^n$. We recall that, because of the algebraic independence of $p_1(x),\ldots,p_n(x)$, there is no polynomial $\widehat f(p)$, not equal to the null polynomial $0$ (that one with all vanishing coefficients), for which $\widehat f(p_1(x),\ldots,p_n(x))=0$, for all $x\in\R^n$.

Given the irreducible finite reflection group $W$, there are infinitely many possible choices of a set of $n$ basic invariant polynomials, but their degrees $d_a={\rm deg}(p_a(x))$, $a=1,\ldots,n$, are well known and typical of the given group $W$. They were determined by Coxeter in Ref.~\cite{cox1951}, and they are all different, except in the case of the groups of type $D_n$, with even $n$, in which case there are two basic invariant polynomials of degree $n$.
Usually, as we do, the basic invariant polynomials are ordered according to their degrees, by requiring
$$ d_1\leq d_2\leq\ldots \leq d_n\,.
$$

The irreducibility of $W$ implies $d_1=2$, and the orthogonality of $W$ allows one to take
\begin{equation}\label{quadraticinv}
  p_1(x)=\sum_{i=1}^n {x_i}^2\,.
\end{equation}
The other $n-1$ basic invariant polynomials have expressions that depend on the specific group, but are not univocally determined. To be more precise, a rank $n$ finite orthogonal reflection group $W$ is generated by the $n$ orthogonal reflections with respect to the hyperplanes orthogonal to the $n$ simple roots. A rotation of the $n$ simple roots with respect to the system of coordinates used in $\R^n$, imply different matrix expressions for the group elements and these imply different explicit expressions of the basic invariant polynomials of degree greater than 2. Moreover, even if the matrix representation of the group $W$ is given, the choice of the basic invariant polynomials is not unique, because it is possible to make algebraic combinations of a given set of basic polynomials to obtain a new set of basic invariant polynomials equally good for that matrix representation (such a change of basis is described in detail in step~\ref{trasformabase} in Section \ref{algorithm}).\\

Several authors determined explicit bases of invariant polynomials for (some of) the irreducible finite reflection groups, for example Coxeter in Ref.~\cite{cox1951}, Ignatenko in Ref.~\cite{Ignatenko1984}, Mehta in Ref.~\cite{Mehta1988}, and many others that are not cited here. A few authors tried to select, among the infinitely many possibilities, some distinguished bases of invariant polynomials, by requiring some supplementary conditions to be satisfied. For example,
Saito, Yano and Sekiguchi in Ref.~\cite{SYS1980} defined the {\em flat bases}, Sartori and Talamini in Ref.~\cite{Sar-Tal1991} the {\em $a$-bases}, and Flatto in Ref.~\cite{Flatto1970} and Iwasaki in Ref.~\cite{Iwasaki1997} the {\em canonical bases}.
In this article we are concerned with the canonical bases, so it is convenient to recall some definitions and results.\\

In Refs.~\cite{Flatto1968} and \cite{FlattoWiener1969} Flatto and Wiener were concerned with certain mean value problems and they established an algorithm to compute bases of invariant polynomials of finite reflection groups (not of type $D_n$, even $n$) 
related to the mean value problems considered. Their algorithm is quite difficult to use for an explicit computation of the bases of invariant polynomials they were interested in. In Ref.~\cite{Flatto1970} Flatto refined his previous results and gave a much simpler algorithm, that can be applied also the groups of type $D_n$, even $n$. To review the algorithm of Ref.~\cite{Flatto1970}, let us define a bilinear map that maps two polynomials $p(x)$ and $q(x)$ into a polynomial $( p,q)(x)$, obtained in the following way:
\begin{equation}\label{Flattoprod}
  ( p,q)(x)=p(\partial)\,q(x)\,,
\end{equation}
where $p(\partial)$ means that in $p(x)$ one has to replace all occurrences of $x_i^k$ with $\frac{\partial^k}{\partial x_i^k}$, $i=1,\ldots,n$, $k\in\mathbb{N}$, and  $p(\partial)\,q(x)$ in Eq. (\ref{Flattoprod}) is the polynomial obtained by the application of the differential operator $p(\partial)$ to the polynomial $q(x)$. If $p(x)$ and $q(x)$ are real homogeneous polynomials of degree $d_p$ and $d_q$, respectively, then $( p,q)(x)$ is 0 if $d_p>d_q$, a real number (also equal to $(q,p)(x)$) if $d_p=d_q$, and a real homogeneous polynomial of degree $d_q-d_p$ if $d_p<d_q$.
In Ref.~\cite{Flatto1970} Flatto proved that, for all irreducible finite reflection groups $W$, the basic invariant polynomials he was interested in, let us call them $h_1(x),\ldots,h_n(x)$, satisfy the following set of partial differential equations: 
\begin{equation}\label{flattocondition1}
  ( h_a,h_b)(x)=0\,,\qquad a,b=1,\ldots,n,\qquad  a\neq b\,.
\end{equation}
For what it is said above, these equations are not trivial only if $d_a\leq d_b$ (with the case $d_a= d_b$ concerning the case $W=D_n$, even $n$, only).
Moreover (see Ref.~\cite{Flatto1970}, pp. 555--556), given a set of basic invariant polynomials $h_1(x),\ldots,h_n(x)$, satisfying Eq. (\ref{flattocondition1}), if $h_a(x)$ is a polynomial with degree $d_a$ of multiplicity one in the set of degrees $\{d_1,\ldots,d_n\}$ (this condition excludes the two basic invariant polynomials of degree $n$ that exist in the case $W=D_n$, even $n$, only), then the polynomial $h_a(x)$ is determined up to real constant multiples. In the case of $D_n$, even $n$, given any two basic invariant polynomials of degree $n$ satisfying Eq. (\ref{flattocondition1}), let us call them $h_{\frac{n}{2}}(x)$ and $h_{\frac{n}{2}+1}(x)$, all other possible couples of invariant polynomials of degree $n$, satisfying Eq. (\ref{flattocondition1}), belong to the linear space $L_n=\{c_1\,h_{\frac{n}{2}}(x)+c_2\,h_{\frac{n}{2}+1}(x)\, \mid\, c_1,c_2\in\R\}$.
Iwasaki in Ref.~\cite{Iwasaki1997} called {\em canonical} a basis of invariant polynomials satisfying Eq. (\ref{flattocondition1}).\\

In Ref.~\cite{Naka-Tsu2014} Nakashima and Tsujie gave a method to determine the canonical bases of the real irreducible finite reflection groups that is different from those proposed in Refs.~\cite{Flatto1968}, \cite{FlattoWiener1969} and \cite{Flatto1970}. Together with Terao, in Ref.~\cite{NTT2016}, they also defined the canonical bases for the unitary reflection groups, and extended to the unitary reflection groups many results of Ref.~\cite{Naka-Tsu2014}.\\

At page 556 of Ref.~\cite{Flatto1970}, Flatto noted that Eq. (\ref{flattocondition1}) allows one to determine by recursion a canonical basis $h_1(x),\ldots,h_n(x)$, going up with the degrees, from $d_1$ to $d_n$. However, this method is not suitable to determine the canonical bases in an explicit way, because Eq. (\ref{flattocondition1}) is non-linear in the unknown invariant polynomials $h_1(x),\ldots,h_n(x)$, and one soon encounters very difficult calculations. In Ref.~\cite{Iwasaki1997}, using the results of Ref.~\cite{Flatto1970}, Iwasaki noted that one can simplify Eq. (\ref{flattocondition1}), by replacing it with the following equivalent set of partial differential equations (equivalent if $W\neq D_n$, even $n$):
\begin{equation}\label{linflattocondition1}
  ( p_a,h_b)(x)=0\,,\qquad a,b=1,\ldots,n,\qquad d_a< d_b\,,
\end{equation}
where $p_1(x),\ldots,p_n(x)$ is any known basis of homogeneous invariant polynomials. In the case of $D_n$, even $n$, a system equivalent to Eq. (\ref{flattocondition1}), is obtained by adding to Eq. (\ref{linflattocondition1}), the partial differential equation:
\begin{equation}\label{linflattocondition2}
    ( h_{\frac{n}{2}},h_{\frac{n}{2}+1})(x)=0\,
    .
\end{equation}
To determine explicitly a canonical basis with this method, one has just to write the most general basis transformation $p\to h$ (see details in step \ref{trasformabase} in Section \ref{algorithm}), and solve Eq. (\ref{linflattocondition1}), and possibly Eq. (\ref{linflattocondition2}), going up with the degrees, from $d_1$ to $d_n$.
The solution of the system of equations arising from Eq. (\ref{linflattocondition1}), and possibly Eq. (\ref{linflattocondition2}), is much simpler than that one arising from Eq. (\ref{flattocondition1}).\\

It is natural to fix a normalization of the polynomials in a canonical basis.
Given a real homogeneous polynomial $p(x)$, let us call  {\em (canonical) norm} of $p$ the non-negative number:
$$\|p\|_c=\sqrt{( p,p)(x)}\,.$$
$\|p\|_c$ is a real positive number if $p(x)\neq 0$, and $0$ if $p(x)= 0$, where $0$ here is the null polynomial.

A normalized canonical basis can be defined, for example, by requiring the following conditions, in addition to those of Eq. (\ref{flattocondition1}) (or to those of Eq. (\ref{linflattocondition1}), and possibly Eq. (\ref{linflattocondition2})):
\begin{equation}\label{norma}
 \|h_a\|_c=1\,,\qquad a=1,\ldots,n\,.
\end{equation}
This normalization was first proposed by Nakashima and Tsujie in Ref.~\cite{Naka-Tsu2014}.
It turns out that a polynomial of a normalized canonical basis is univocally determined in all cases in which its degree has multiplicity one in the set of degrees $\{d_1,\ldots,d_n\}$, and that, in the case of groups of type $D_n$, even $n$, the two equal degree polynomials of a normalized canonical basis are determined up to orthogonal transformations in the linear space $L_n$ they generate.\\

The explicit expressions of the polynomials in a normalized canonical basis have real coefficients that usually are irrational and have non-trivial common divisors. In addition, $p_1(x)$ has not the standard form (\ref{quadraticinv}). For these reasons, in some applications it might be better to avoid the normalization of Eq. (\ref{norma}).\\

Because of the special form of $p_1(x)$ (that is given by Eq. (\ref{quadraticinv}), or is proportional to it, as in the normalized canonical bases), the differential operator $p_1(\partial)$ is the Laplace operator $\Delta$ (or is proportional to it), in fact:
$$p_1(x)=\sum_{i=1}^n\,x_i^2\qquad \Rightarrow \qquad p_1(\partial)=\sum_{i=1}^n\,\frac{\partial^2}{\partial x_i^2}=\Delta$$
This implies that all basic polynomials $h_b(x)$, $b=2,\ldots,n$, satisfying Eq. (\ref{flattocondition1}), satisfy in particular the equations $( h_1,h_b)(x)=0$, $b=2,\ldots,n$, that is, the Laplace equations $$\Delta\,h_b(x)=0\,,\qquad b=2,\ldots,n\,,$$ and are thus harmonic polynomials.
This means, in particular, that the basic polynomials $h_b(x)$, $b=2,\ldots,n$, satisfy all nice theorems on harmonic functions, like the mean value property, the maximum-minimum principle and so on.\\

The explicit form of the canonical bases of the irreducible finite reflection groups of types $A_n$, $B_n$, $D_n$ and $I_2(m)$, were determined in Ref.~\cite{Iwasaki1997} and those of type $H_3$, $H_4$ and $F_4$ were determined in Ref.~\cite{IwasakiKenmaMatsumoto2002}. The explicit form of the canonical bases of the irreducible finite reflection groups of types $E_6$, $E_7$ and $E_8$ are  determined in this article. The author verified the results of Ref.~\cite{IwasakiKenmaMatsumoto2002} and found that they are correct, but for the sign of the monomials containing $\Delta_3$ in Table 4, concerning the basic invariant polynomials of $H_3$. The canonical bases in Ref.~\cite{IwasakiKenmaMatsumoto2002} are not normalized and to get the normalized canonical bases one has to multiply each of the basic polynomials there listed by a suitable real factor.\\

\section{The algorithm\label{algorithm}}

The algorithm the author used to determine the explicit canonical bases for the groups of types $E_6$, $E_7$ and $E_8$, is essentially the same as that one used in Ref.~\cite{IwasakiKenmaMatsumoto2002}. It is described in detail in the following 6 steps. The algorithm here presented can be used to find the explicit canonical bases of any irreducible finite reflection groups not of type $D_n$, even $n$. 
At the end of this section it will be described how one can modify the algorithm to apply it to the groups of type $D_n$, even $n$.\\

\begin{enumerate}
\item One starts with a known basis, $p_a(x)$, $a=1,\ldots,n$, of $W$-invariant homogeneous polynomials ($W$ is here any irreducible finite reflection groups not of type $D_n$, even $n$). One can use, for example, the bases determined by Mehta in Ref.~\cite{Mehta1988}. These bases are not canonical.\\

\item \label{trasformabase} One writes the most general basis transformation $p_a(x)\to h_a(x)$, $a=1,\ldots,n$, to get a new basis of $W$-invariant polynomials $h_a(x)$, $a=1,\ldots,n$. In order to obtain the new polynomials $h_a(x)$,  $a=1,\ldots,n$, forming a set of algebraically independent homogeneous polynomials of degrees $d_a$, $a=1,\ldots,n$, there are some conditions to satisfy. As the new basic polynomials $h_a(x)$,  $a=1,\ldots,n$, are real $W$-invariant polynomials, it must be possible to write them as real
polynomials of the original basic invariant polynomials $p_a(x)$, $a=1,\ldots,n$. One has then expressions like the following ones:
$$  h_a(x)={\widehat h}_a(p_1(x),\ldots,p_n(x)),\qquad a=1,\ldots,n\,,
$$
where the ${\widehat h}_a(p_1,\ldots,p_n)={\widehat h}_a(p)$ are real polynomial functions of $p\in\R^n$. 
 The homogeneity of the $h_a(x)$ and the set of degrees $d_1,\ldots,d_n$, of $W$, restrict the possible monomials appearing in ${\widehat h}_a(p)$: a monomial $p_1^{m_1}\,p_2^{m_2}\cdots p_n^{m_n}$ is in ${\widehat h}_a(p)$ if and only if $\sum_{b=1}^n d_b\,m_b=d_a$. So one can easily write down explicitly the most general expressions of the ${\widehat h}_a(p)$, $a=1,\ldots,n$. Let us call $z_i$, $i=1,2,\ldots,$ the (unknown) real coefficients of the monomials in ${\widehat h}_a(p)$. For example, for $W=E_6$, the degrees of the basic invariant polynomials are $d_1,\ldots,d_6=2,5,6,8,9,12$, and one writes ${\widehat h}_6(p)=z_1\, p_1^6 + z_2 \,p_1 p_2^2 + z_3\, p_1^3 p_3 + z_4\, p_3^2 + z_5 \,p_1^2 p_4 + z_6\, p_6$. If one requires that the expression of ${\widehat h}_a(p)$ necessarily contains the variable $p_a$, for each $a=1,\ldots,n$, one is sure that the transformed invariant polynomials are algebraically independent. This also implies that the inverse transformation $h\to p$ is well defined and polynomial.\\

\item One substitutes in ${\widehat h}_a(p)$ the variables $p_b$, $b=1,\ldots,n$, with the basic invariant polynomials $p_b(x)$, $b=1,\ldots,n$, obtaining the real polynomial functions $h_a(x)$. The coefficients of the monomials in $x_1,\ldots,x_n$ that appear in $h_a(x)$, are real linear combinations of the $z_i$'s.\\

\item One inserts these general expressions of the polynomials $h_a(x)$, $a=1,\ldots,n$, in Eq. (\ref{linflattocondition1}). The coefficients of similar monomials (with respect to the variables $x_1,\ldots,x_n$) give the equations of a (quite big) system of linear equations for the indeterminates $z_i$'s.\\

\item One solves this linear system, and obtains the values of the $z_i$'s that define the basis transformation from the original basis $p_a(x)$, $a=1,\ldots,n$, to the canonical basis $h_a(x)$, $a=1,\ldots,n$. In each of the $n$ explicit expressions $\widehat h_a(p)$, $a=1,\ldots,n$, one of the indeterminates $z_i$ remains undetermined and appears as a common multiplicative factor. If it is not required to find the normalized canonical basis, one can give to these $n$ indeterminates $z_i$ some convenient real values in such a way to avoid denominators and common divisors in the explicit expressions of the $h_a(x)$, $a=1,\ldots,n$.\\

\item To determine the normalized canonical basis $k_1(x),\ldots,k_n(x)$, from a non-normalized canonical basis $h_1(x),\ldots,h_n(x)$, one first defines $k_a(x)=z_a\,h_a(x)$, $a=1,\ldots, n$, with $z_1,\ldots,z_n$, arbitrary real coefficients, and then solves the non-linear system in the $z_a$'s obtained by substituting the polynomials $k_a(x)$, $a=1,\ldots,n$, in Eq. (\ref{norma}), that is, by solving the system $\ \|k_a\|_c=1$, $\ a=1,\ldots,n$.\\
    \end{enumerate}

In the case of groups of type $D_n$, even $n$, one can use the same algorithm presented above, with the following slight modifications at steps 4, 5 and 6. Step 4: In addition to Eq. (\ref{linflattocondition1}), one has to consider also Eq.~(\ref{linflattocondition2}) that yields a system of non-linear equations in the indeterminates $z_i$'s. Step 5: Instead of $n$, there are $n+1$ variables $z_i$, that remain undetermined. 
If one is not interested to find the normalized canonical basis, one can give to these $n+1$ variables any convenient real values, otherwise, one can normalize the polynomials so found using Eq. (\ref{norma}), as written in step 6. Step 6: If one has to determine a normalized canonical basis from a non-normalized canonical basis $h_a(x)$, $a=1,\ldots,n$, one defines
$k_a(x)$, $a=1,\ldots,n$, $a\neq {\frac{n}{2}}, {\frac{n}{2}+1}$, as written in step 6, and defines $k_{\frac{n}{2}}(x)=c_{1}\,h_{\frac{n}{2}}(x)+c_{2}\,h_{\frac{n}{2}+1}(x)$ and $k_{\frac{n}{2}+1}(x)=c_{3}\,h_{\frac{n}{2}}(x)+c_{4}\,h_{\frac{n}{2}+1}(x)$, with $c_{1},c_2,c_3,c_{4}$ arbitrary real coefficients. 
One then inserts $k_{\frac{n}{2}}(x)$ and $k_{\frac{n}{2}+1}(x)$ in Eq. (\ref{linflattocondition2}), determining one of the $c_i$'s, and normalizes all $k_a(x)$, $a=1,\ldots,n$, as explained in step 6. After all this, one of the variables $c_{1},c_2,c_3,c_{4}$, remains still undetermined, and one can give it any convenient real value.\\

To get concretely the results the author used the commercial computer algebra system {\em Mathematica}.\\

The explicit form of the canonical bases, as well as those of the normalized canonical bases, for the irreducible finite reflection groups of type $E_6$, $E_7$ and $E_8$, are specified in the next three sections, respectively.
As the explicit expressions of the canonical basic polynomials $h_a(x)$, $a=1,\ldots,n$, are very large, it will be specified just the basis transformations to get them from a known basis of invariant polynomials, and not written their explicit expressions in full.

\section{Canonical basis of invariant polynomials for $E_6$\label{E6}}
In Ref.~\cite{Mehta1988} Mehta described how to obtain a basis of invariant polynomials for $E_6$. Consider $x=(x_1,\ldots,x_6)$ and the 27 linear forms ${\ell}_{1}(x),\ldots,\ell_{27}(x)$, listed in (4.2) of Ref.~\cite{Mehta1988}, copied here under:
$$2\sqrt{\frac{2}{3}}\,x_6\,,\qquad \sqrt{\frac{2}{3}}\left(\pm\sqrt{3}\,x_5-x_6\right)\,,$$
$$\pm x_3\pm x_4-\sqrt{\frac{2}{3}}\,x_6\,,\qquad \pm x_1\pm x_2-\sqrt{\frac{2}{3}}\,x_6\,,$$
$$\pm x_2\pm x_4+\frac{1}{\sqrt{6}}\left(\sqrt{3}\,x_5+x_6\right)\,,\qquad \pm x_1\pm x_3+\frac{1}{\sqrt{6}}\left(\sqrt{3}\,x_5+x_6\right)\,,$$
$$\pm x_2\pm x_3-\frac{1}{\sqrt{6}}\left(\sqrt{3}\,x_5-x_6\right)\,,\qquad \pm x_1\pm x_4-\frac{1}{\sqrt{6}}\left(\sqrt{3}\,x_5-x_6\right)\,.$$
A basis of algebraically independent homogeneous polynomials for the orthogonal irreducible group of type $E_6$ is obtained in the following way:
\begin{equation}\label{basisE6}
  p_a(x)=\sum_{k=1}^{27}\,(\ell_k(x))^{d_a}\,,\qquad a=1,\ldots,6\,,
\end{equation}
where $(d_1,\ldots,d_6)=(2,\,5,\,6,\,8,\,9,\,12)$ are the degrees of the basic invariant polynomials of $E_6$.\\
From the basic invariant polynomials (\ref{basisE6}), it is possible to calculate the 36 reflecting hyperplanes, that are also given in Eq. (4.1) of Ref.~\cite{Mehta1988}, whose gradients specify the root directions. The method to do this derives from a classical result of Coxeter, Theorem 6.2 of Ref.~\cite{cox1951}. One first calculates the jacobian matrix $j(x)$, whose elements are $j_{ai}(x)=\frac{\partial p_a(x)}{\partial x_i}$, $a,i=1,\ldots,n$. The determinant $\det(j(x))$ is proportional to the product of the $N=\sum_{a=1}^n (d_a-1)$ linear forms $L_r(x)$, $r=1,\ldots,N$, which determine the equations $L_r(x)=0$ of the reflecting hyperplanes, and whose gradients $\nabla L_r(x)$ are proportional to the roots. In Ref.~\cite{Mehta1988} Mehta did not use this method to find the equations of the reflecting hyperplanes but instead he calculated both the equations of the reflecting hyperplanes and the explicit expressions of the basic invariant polynomials from the linear forms ${\ell}_{1}(x),\ldots,\ell_{27}(x)$, once he proved that they form an invariant set. A possible set of simple roots is the following one:
%
$$\alpha_1= e_2- e_3+\frac{1}{\sqrt{2}}\, e_5+\sqrt{\frac{3}{2}}\, e_6,\qquad \alpha_2= e_3- e_4-\sqrt{2}\, e_5,\qquad
\alpha_3=2\, e_4,$$
$$\alpha_4= e_3- e_4+\sqrt{2}\, e_5,\qquad \alpha_5= e_2- e_3-\frac{1}{\sqrt{2}}\, e_5-\sqrt{\frac{3}{2}}\, e_6,\qquad \alpha_6= e_1- e_2- e_3- e_4,
$$
where the $ e_i$, $i=1,\ldots,6$, are the canonical unit vectors of $\R^6$. The simple roots just specified take place in the following Coxeter-Dynkin diagram:
   \begin{center}
 \unitlength=1.0mm
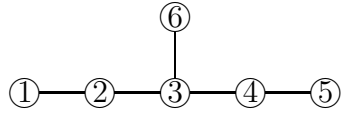
 \begin{picture}(40,15)(25,-3)
    \put(20,0){\circle{4}}
    \put(30,0){\circle{4}}
    \put(40,0){\circle{4}}
    \put(50,0){\circle{4}}
    \put(60,0){\circle{4}}
    \put(40,10){\circle{4}}


    \put(22,0){\line(1,0){6}}
    \put(32,0){\line(1,0){6}}
    \put(42,0){\line(1,0){6}}
    \put(52,0){\line(1,0){6}}
    \put(40,2){\line(0,1){6}}

 \put(19.0,-1.2){$1$}
 \put(29.0,-1.2){$2$}
 \put(39.0,-1.2){$3$}
 \put(49.0,-1.2){$4$}
 \put(59.0,-1.2){$5$}
 \put(39.0,8.8){$6$}

\end{picture}
    \end{center}
where the numbers in the circles correspond to the simple root indices. The reflections corresponding to the simple roots $\alpha_1,\ldots,\alpha_6$ generate the group of type $E_6$ for which the polynomials in Eq. (\ref{basisE6}) are invariant.\\

The basic invariant polynomials given by Eq. (\ref{basisE6}) may not be the best ones for computations, as they contain rational and irrational numbers for coefficients, and $p_1(x)$ is not as in Eq. (\ref{quadraticinv}), so one might prefer to consider some multiples of them, at least to get rid of the numeric common divisors and denominators. A possible choice is to take as basic invariant polynomials the following ones (the argument $x$ of the polynomials $p_a(x)$, $q_a(x)$, $a=1,\ldots,n$, and of all polynomials appearing in the rest of this section, is understood):
$$q_1=\frac{1}{12}\,p_1 ,\ \ q_2 =\frac{3}{20}\sqrt{2}\,p_2 ,\ \ q_3 =\frac{1}{2}\,p_3 ,\ \ q_4 =9\,p_4 ,\ \ q_5 =\frac{9}{7}\sqrt{2}\,p_5 ,\ \ q_6 =36\,p_6.$$

A basis transformation to get a canonical basis $h_1,\ldots,h_6$ is the following one:
{\footnotesize
$$h_1=q_1\,,\quad h_2=\frac{1}{\sqrt{3}}\,q_2\,,\quad h_3=\sqrt{3}\,(-8\,q_1^3+q_3)\,,$$ $$h_4=\frac{\sqrt{3}}{5}\,(1120\,q_1^4-224\,q_1q_3+3\,q_4)\,,\quad h_5=\frac{1}{\sqrt{3}}\,(-80\,q_1^2q_2+q_5)\,,$$
$$h_6=\frac{1}{405}\,(-169845984\,q_1^6 - 18714080\,q_1q_2^2 +
  50516928\,q_1^3q_3 - 657888\,q_3^2 -
  1108536\,q_1^2 q_4 + 21171\,q_6)\,.$$}

The normalized canonical basis $k_1,\ldots,k_6$, satisfying Eq. (\ref{norma}), is obtained with the following multiples of the polynomials 
$h_1,\ldots,h_6$:
$$k_1=\frac{1}{2\,\sqrt{3}}\,h_1\,,\qquad k_2=\frac{1}{48\,\sqrt{2}}\,h_2\,,\qquad k_3=\frac{1}{576\,\sqrt{5}}\,h_3\,,$$
$$k_4=\frac{1}{13824\,\sqrt{70}}\,h_4\,,\qquad k_5=\frac{1}{46080\,\sqrt{2}}\,h_5\,,\qquad k_6=\frac{1}{4423680\,\sqrt{543389}}\,h_6\,.$$

\section{Canonical basis of invariant polynomials for $E_7$}
In Ref.~\cite{Mehta1988} Mehta described how to obtain a basis of invariant polynomials for $E_7$. Consider $x=(x_1,\ldots,x_7)$ and the 56 linear forms $\ell_{1}(x),\ldots,\ell_{56}(x)$, listed in (4.4) of Ref.~\cite{Mehta1988}, copied here under:
$$\pm x_i\pm x_j\pm x_k\,,$$
where the indices $(i,j,k)$ can take the following $7$ sets of values: $(i,j,k)=(1,2,7),$ $(1,3,6),$ $(1,4,5),$ $(2,3,5),$ $(2,4,6),$ $(3,4,7),$ $(5,6,7)$.

A basis of algebraically independent homogeneous polynomials for the orthogonal irreducible group of type $E_7$ is obtained in the following way:
\begin{equation}\label{basisE7}
  p_a(x)=\sum_{k=1}^{56}\,(\ell_k(x))^{d_a}\,,\qquad a=1,\ldots,7\,,
\end{equation}
where $(d_1,\ldots,d_7)=(2,\,6,\,8,\,10,\,12,\,14,\,18)$ are the degrees of the basic invariant polynomials of $E_7$. (All the degrees are even, so it would be sufficient in Eq. (\ref{basisE7}) to sum over the 28 linear forms that start with the plus sign).\\
From the basic invariant polynomials (\ref{basisE7}), 
it is possible to calculate the 63 reflecting hyperplanes, that are also given in Eq. (4.3) of Ref.~\cite{Mehta1988}, whose gradients specify the root directions. (The method to do this was described in Section \ref{E6}). A possible set of simple roots is the following one:
$$\alpha_1=2\, e_7,\qquad \alpha_2= e_2- e_3- e_6- e_7,\qquad
\alpha_3=2\, e_6,$$
$$\alpha_4= e_3- e_4- e_5- e_6,\qquad \alpha_5=2\, e_4,\qquad \alpha_6= e_1- e_2- e_3- e_4,\qquad \alpha_7=2\, e_5,
$$
where the $ e_i$, $i=1,\ldots,7$, are the canonical unit vectors of $\R^7$. The simple roots just specified take place in the following Coxeter-Dynkin diagram:
   \begin{center}
 \unitlength=1.0mm
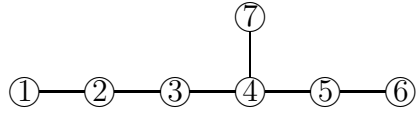
 \begin{picture}(40,15)(23,-3)
    \put(10,0){\circle{4}}
    \put(20,0){\circle{4}}
    \put(30,0){\circle{4}}
    \put(40,0){\circle{4}}
    \put(50,0){\circle{4}}
    \put(60,0){\circle{4}}
    \put(40,10){\circle{4}}


    \put(12,0){\line(1,0){6}}
    \put(22,0){\line(1,0){6}}
    \put(32,0){\line(1,0){6}}
    \put(42,0){\line(1,0){6}}
    \put(52,0){\line(1,0){6}}
    \put(40,2){\line(0,1){6}}

  \put(9.0,-1.2){$1$}
 \put(19.0,-1.2){$2$}
 \put(29.0,-1.2){$3$}
 \put(39.0,-1.2){$4$}
 \put(49.0,-1.2){$5$}
 \put(59.0,-1.2){$6$}
 \put(39.0,8.8){$7$}

\end{picture}
    \end{center}
where the numbers in the circles correspond to the simple root indices.
The reflections corresponding to the simple roots $\alpha_1,\ldots,\alpha_7$ generate the group of type $E_7$ for which the polynomials in Eq. (\ref{basisE7}) are invariant.\\


We start from the basic invariant polynomials given by the following multiples of those in Eq. (\ref{basisE7}):
$$q_i=\frac{1}{24}\,p_i,\ \ i=1,2,4,5,7,\qquad q_j=\frac{1}{8}\,p_j,\ \ j=3,6.$$
The explicit expressions of the basic invariant polynomials $q_1(x),\ldots,q_7(x)$, together with the computer code used to get them, are available online in the supplemental material of Ref.~\cite{tal-jmp2010}.\\

A basis transformation to get a canonical basis $h_1,\ldots,h_7$ is the following one:
{\footnotesize
$$h_1=q_1\,,\qquad h_2=\frac{1}{2}\,(-15 \,q_1^3+11\,q_2)\,,\qquad h_3=\frac{1}{20}\,(2835 \,q_1^4 - 3276 \,q_1 q_2 + 247 \,q_3)\,,$$
$$h_4=\frac{1}{2}\,\big(-9 \,q_1 \,(18 \,q_1^4 - 30 \,q_1 q_2 + 5 \,q_3) + 23 \,q_4\big)\,,$$ $$h_5=\frac{1}{10}\,\big(11\,(280 \,q_1^6 - 1288 \,q_1^3 q_2 - 490 \,q_2^2 +
       761 \,q_1^2 q_3 - 970 \,q_1 q_4) + 1735 \,q_5\big)\,,$$
$$h_6=\frac{1}{319}\,\Big(819\,\big(33\,(4490\,q_1^7 - 8666\,q_1^4 q_2 + 4900 \,q_1 q_2^2 -
         300 \,q_1^3 q_3 - 465 \,q_2 q_3 + 2525 \,q_1^2 q_4) \,+$$
         $$-\,
       36115 \,q_1 q_5\big) + 1610605 \,q_6\Big)\,,$$
$$h_7=\frac{1}{15682040}\,\Big(-2431\,\big(5085078551185 \,q_1^9 - 11402026037640 \,q_1^6 q_2 +
       7472423123536 \,q_1^3 q_2^2 \,+$$
       $$-\, 201739938400 \,q_2^3 -
       540102070990 \,q_1^5 q_3- 748116822184 \,q_1^2 q_2 q_3\,+$$
       $$  -\,
       46311340011 \,q_1 q_3^2 + 40 \,(152224768729 \,q_1^4 -
         6555491354 \,q_1 q_2 + 1476892164 \,q_3) \,q_4\big)\, +$$
$$+\,     742560\,(12033352910 \,q_1^3 - 517191829 \,q_2) \,q_5 -
     703975240263600 \,q_1^2 q_6 + 64758924763060 \,q_7\Big)\,.$$}

     The normalized canonical basis $k_1,\ldots,k_7$, satisfying Eq. (\ref{norma}), is obtained with the following multiples of the polynomials 
$h_1,\ldots,h_7$:
$$k_1=\frac{1}{\sqrt{14}}\,h_1\,,\qquad k_2=\frac{1}{24\,\sqrt{2310}}\,h_2\,,\qquad k_3=\frac{1}{2016\,\sqrt{741}}\,h_3\,,$$
$$k_4=\frac{1}{40320\,\sqrt{138}}\,h_4\,,\qquad k_5=\frac{1}{483840\,\sqrt{7634}}\,h_5\,,$$ $$k_6=\frac{1}{11612160\,\sqrt{146565055}}\,h_6\qquad k_7=\frac{1}{92897280\,\sqrt{5181830514230370}}\,h_7\,.$$

\section{Canonical basis of invariant polynomials for $E_8$}
In Ref.~\cite{Mehta1988} Mehta described how to obtain a basis of invariant polynomials for $E_8$. Consider $x=(x_1,\ldots,x_8)$ and the 240 linear forms $\ell_{1}(x),\ldots,\ell_{240}(x)$, listed at page 1097 of Ref.~\cite{Mehta1988}, copied here under:
$$\pm 2\,x_i\,,\qquad i=1,\ldots,8\,,$$
$$\pm x_i\pm x_j\pm x_k\pm x_l\,,$$
where the indices $(i,j,k,l)$ can take the following 14 sets of values: $(i,j,k,l)=(1,2,3,4),$ $(1,2,5,6),$ $(1,2,7,8),$ $(1,3,5,7),$ $(1,3,6,8),$ $(1,4,6,7),$ $(1,4,5,8),$ $(2,3,5,8),$ $(2,3,6,7),$ $(2,4,5,7),$ $
(2,4,6,8),$ $(3,4,5,6),$ $(3,4,7,8),$ $(5,6,7,8)$.
A basis of algebraically independent homogeneous polynomials for the orthogonal irreducible group of type $E_8$ is obtained in the following way:
\begin{equation}\label{basisE8}
  p_a(x)=\sum_{k=1}^{240}\,(\ell_k(x))^{d_a}\,,\qquad a=1,\ldots,8\,,
\end{equation}
where $(d_1,\ldots,d_8)=(2,\,8,\,12,\,14,\,18,\,20,\,24,\,30)$ are the degrees of the basic invariant polynomials of $E_8$.  (All the degrees are even, so it would be sufficient in Eq. (\ref{basisE8}) to sum over the 120 linear forms that start with the plus sign).\\
From the basic invariant polynomials (\ref{basisE8}), 
it is possible to calculate the 120 reflecting hyperplanes, whose equations coincide with $\ell_{j}(x)=0$, $j=1,\ldots, 240$, (counting twice each equation), whose gradients specify the root directions. (The method to do this was described in Section \ref{E6}). A possible set of simple roots is the following one:
$$\alpha_1= e_1- e_2- e_3- e_4,\qquad \alpha_2=2\, e_4,\qquad
\alpha_3= e_3- e_4- e_5- e_6,$$
$$\alpha_4=2\, e_6,\qquad \alpha_5= e_5- e_6- e_7- e_8,\qquad \alpha_6=2\, e_8,\qquad \alpha_7= e_2- e_3- e_5- e_8,$$ $$\alpha_8=2\, e_7,
$$
where the $ e_i$, $i=1,\ldots,8$, are the canonical unit vectors of $\R^8$. The simple roots just specified take place in the following Coxeter-Dynkin diagram:
    \begin{center}
 \unitlength=1.0mm
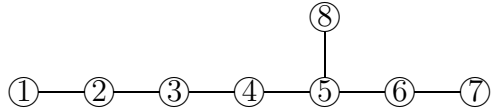
 \begin{picture}(40,15)(20,-3)
    \put(0,0){\circle{4}}
    \put(10,0){\circle{4}}
    \put(20,0){\circle{4}}
    \put(30,0){\circle{4}}
    \put(40,0){\circle{4}}
    \put(50,0){\circle{4}}
    \put(60,0){\circle{4}}
    \put(40,10){\circle{4}}


    \put(2,0){\line(1,0){6}}
    \put(12,0){\line(1,0){6}}
    \put(22,0){\line(1,0){6}}
    \put(32,0){\line(1,0){6}}
    \put(42,0){\line(1,0){6}}
    \put(52,0){\line(1,0){6}}
    \put(40,2){\line(0,1){6}}

 \put(-1,-1.2){$1$}
 \put(9.0,-1.2){$2$}
 \put(19.0,-1.2){$3$}
 \put(29.0,-1.2){$4$}
 \put(39.0,-1.2){$5$}
 \put(49.0,-1.2){$6$}
 \put(59.0,-1.2){$7$}
 \put(39.0,8.8){$8$}

\end{picture}
    \end{center}
where the numbers in the circles correspond to the simple root indices. The reflections corresponding to the simple roots $\alpha_1,\ldots,\alpha_8$ generate the group of type $E_8$ for which the polynomials in Eq. (\ref{basisE8}) are invariant.\\

We start from the basic invariant polynomials given by the following multiples of those in Eq. (\ref{basisE8}):
$$q_1=\frac{1}{120}\,p_1,\qquad q_i=\frac{1}{48}\,p_i,\ \ i=2,3,4,5,6,7,\qquad q_8=\frac{1}{240}\,p_8.$$
The explicit expressions of the basic invariant polynomials $q_1(x),\ldots,q_8(x)$, together with the computer code used to get them, are available online in the supplemental material of Ref.~\cite{tal-jmp2010}.\\

A basis transformation to get a canonical basis $h_1(x),\ldots,h_8(x)$ is the following one:

{\footnotesize
$$h_1=q_1\,,\qquad h_2=-10 \,q_1^4+q_2\,,\qquad h_3=\frac{1}{7}\,(4235\,q_1^6 - 495q_1^2\,q_2 + 13\,q_3)\,,$$
$$h_4=\frac{1}{11}\,(-17589\,q_1^7 + 2145\,q_1^3 q_2 -
     91\,q_1 q_3 + 8\,q_4)\,,$$
$$h_5=\frac{1}{7280}\,\big(17\,q_1\,(27922895\,q_1^8 - 3333330\,q_1^4 q_2 -
       24453\,q_2^2 + 227864\,q_1^2 q_3 - 36144\,q_1 q_4) + 7600\,q_5\big)\,,$$
$$h_6=\frac{1}{748}\,\Big(-969\,\big(429\,q_1^2\,(992005\,q_1^8 - 115710\,q_1^4 q_2 -
         1271\,q_2^2)\, +$$
         $$+\, 728\,(5059\,q_1^4 + 10\,q_2)\,q_3 -
       647312\,q_1^3 q_4\big) - 12549880\,q_1 q_5 + 880796\,q_6\Big)\,,$$
$$h_7=\frac{1}{28647880800}\,\bigg(23\,\Big(4199\,\big(20274537662080415\,q_1^{12} - 2250467375658810\,q_1^8 q_2 \,+$$
$$-\, 40769297380581\,q_1^4 q_2^2  + 178341143921528\,q_1^6 q_3+1420510398720\,q_1^2 q_2 q_3\,+$$
$$ -\,
         1640\,(8061383\,q_2^3 + 3080560\,q_3^2)\big) -
       46512\,q_1\,(3015480163976\,q_1^4\, +$$
       $$+\,  10722788425\,q_2)\,q_4+
       3837192062311440\,q_1^3 q_5 - 443023026566400\,q_1^2 q_6\Big) +
     78547609202400\,q_7\bigg)\,,$$
$$h_8=\frac{1}{19626789759713136000000}\bigg(667\,\Big(323\,\big(-11\,q_1^3\,(1094671830559801212459572195245\,q_1^{12}\,+$$
$$ -\,
           124907605937936839186287677130\,q_1^8 q_2 -
           1777752453446126054618835543\,q_1^4 q_2^2 \,+$$
$$+\,           530575867656892216179020\,q_2^3) - 1144\,q_1\,
          (92073813834207882297946571\,q_1^8 \,+$$
           $$+\,511036172390511143554680\,q_1^4 q_2 - 449662651462146636150\,q_2^2)\,q_3 +
         40125576319460176480000\,q_1^3 q_3^2\, +$$
         $$+\,
         48\,(406458546004454701971207148\,q_1^8 + 1088335288483913057060495\,q_1^4 q_2 - 469314879303830560000\,q_2^2\, +$$
         $$+\,
           9038224261602298802000\,q_1^2 q_3)\,q_4 -
         38325110052347105500800\,q_1 q_4^2\big)\,+$$
         $$ -\,
       117040\,(1497805112196088699388741\,q_1^6 - 1376419171002564521550\,q_1^2 q_2 + 9544045308237440000\,q_3)\,q_5 \,+$$
       $$+\,
       100900800\,q_1\,(207536086797307841747\,q_1^4 - 259054452883944920\,q_2)\,q_6\Big) \,+$$
       $$-\, 121281139650829621358815920000\,q_1^3 q_7 +     486032507227341717350400000\,q_8\bigg)\,.$$}

The normalized canonical basis $k_1,\ldots,k_8$, satisfying Eq. (\ref{norma}), is obtained with the following multiples of the polynomials 
$h_1,\ldots,h_8$:
$$k_1=\frac{1}{4}\,h_1\,,\qquad k_2=\frac{1}{1920\,\sqrt{42}}\,h_2\,,\qquad k_3=\frac{1}{92160\,\sqrt{15015}}\,h_3\,,\qquad
k_4=\frac{1}{15482880\,\sqrt{65}}\,h_4\,,$$
$$k_5=\frac{1}{1857945600\,\sqrt{17765}}\,h_5\,,\qquad k_6=\frac{1}{52022476800\,\sqrt{4778475585}}\,h_6\,,$$ $$k_7=\frac{1}{41736889958400\,\sqrt{342348352885}}\,h_7\,,$$
$$k_8=\frac{1}{14383174385664000\,\sqrt{14557753942206761}}\,h_8\,.$$



\end{document}